\documentclass[11pt]{amsart}
\usepackage{amsmath, amsxtra, amsthm, amssymb,eucal,dsfont}
\usepackage[all]{xy}
\usepackage{color}
\usepackage[colorlinks=true, pdfstartview=FitV, linkcolor=blue, citecolor=blue, urlcolor=blue]{hyperref}

\numberwithin{equation}{section}
\swapnumbers

\newtheorem{thm}{Theorem}[section]
\newtheorem{prop}[thm]{Proposition}
\newtheorem{lem}[thm]{Lemma}
\newtheorem{cor}[thm]{Corollary}
\newtheorem{rmk}[thm]{Remark}

\newtheorem{example}[thm]{Example}

\newcommand{\nc}{\newcommand}
\nc{\rc}{\renewcommand}

\rc{\b}{\mathbb}
\rc{\c}{\mathcal}
\nc{\on}{\operatorname}
\nc{\tn}{\textnormal}

\nc{\bA}{\b A}
\nc{\bB}{\b B}
\nc{\bC}{\b C}
\nc{\bD}{\b D}
\nc{\bE}{\b E}
\nc{\bF}{\b F}
\nc{\bG}{\b G}
\nc{\bH}{\b H}
\nc{\bI}{\b I}
\nc{\bJ}{\b J}
\nc{\bK}{\b K}
\nc{\bL}{\b L}
\nc{\bM}{\b M}
\nc{\bN}{\b N}
\nc{\bO}{\b O}
\nc{\bP}{\b P}
\nc{\bQ}{\b Q}
\nc{\bR}{\b R}
\nc{\bS}{\b S}
\nc{\bT}{\b T}
\nc{\bU}{\b U}
\nc{\bV}{\b V}
\nc{\bW}{\b W}
\nc{\bX}{\b X}
\nc{\bY}{\b Y}
\nc{\bZ}{\b Z}

\nc{\cA}{\c A}
\nc{\cB}{\c B}
\nc{\cC}{\c C}
\nc{\cD}{\c D}
\nc{\cE}{\c E}
\nc{\cF}{\c F}
\nc{\cG}{\c G}
\nc{\cH}{\c H}
\nc{\cI}{\c I}
\nc{\cJ}{\c J}
\nc{\cK}{\c K}
\nc{\cL}{\c L}
\nc{\cM}{\c M}
\nc{\cN}{\c N}
\nc{\cO}{\c O}
\nc{\cP}{\c P}
\nc{\cQ}{\c Q}
\nc{\cR}{\c R}
\nc{\cS}{\c S}
\nc{\cT}{\c T}
\nc{\cU}{\c U}
\nc{\cV}{\c V}
\nc{\cW}{\c W}
\nc{\cX}{\c X}
\nc{\cY}{\c Y}
\nc{\cZ}{\c Z}

\nc{\fA}{{\mathfrak A}}
\nc{\fB}{{\mathfrak B}}
\nc{\fC}{{\mathfrak C}}
\nc{\fD}{{\mathfrak D}}
\nc{\fE}{{\mathfrak E}}
\nc{\fF}{{\mathfrak F}}
\nc{\fG}{{\mathfrak G}}
\nc{\fH}{{\mathfrak H}}
\nc{\fI}{{\mathfrak I}}
\nc{\fJ}{{\mathfrak J}}
\nc{\fK}{{\mathfrak K}}
\nc{\fL}{{\mathfrak L}}
\nc{\fM}{{\mathfrak M}}
\nc{\fN}{{\mathfrak N}}
\nc{\fO}{{\mathfrak O}}
\nc{\fP}{{\mathfrak P}}
\nc{\fQ}{{\mathfrak Q}}
\nc{\fR}{{\mathfrak R}}
\nc{\fS}{{\mathfrak S}}
\nc{\fT}{{\mathfrak T}}
\nc{\fU}{{\mathfrak U}}
\nc{\fV}{{\mathfrak V}}
\nc{\fW}{{\mathfrak W}}
\nc{\fZ}{{\mathfrak Z}}
\nc{\fX}{{\mathfrak X}}
\nc{\fY}{{\mathfrak Y}}
\nc{\fa}{{\mathfrak a}}
\nc{\fb}{{\mathfrak b}}
\nc{\fc}{{\mathfrak c}}
\nc{\fd}{{\mathfrak d}}
\nc{\fe}{{\mathfrak e}}
\nc{\ff}{{\mathfrak f}}
\nc{\fg}{{\mathfrak g}}
\nc{\fh}{{\mathfrak h}}
\nc{\fiI}{{\mathfrak i}}  
\nc{\ffi}{{\mathfrak i}}  
\nc{\fj}{{\mathfrak j}}
\nc{\fk}{{\mathfrak k}}
\nc{\fl}{{\mathfrak{l}}}
\nc{\fm}{{\mathfrak m}}
\nc{\fn}{{\mathfrak n}}
\nc{\fo}{{\mathfrak o}}
\nc{\fp}{{\mathfrak p}}
\nc{\fq}{{\mathfrak q}}
\nc{\fr}{{\mathfrak r}}
\nc{\fs}{{\mathfrak s}}
\nc{\ft}{{\mathfrak t}}
\nc{\fu}{{\mathfrak u}}
\nc{\fv}{{\mathfrak v}}
\nc{\fw}{{\mathfrak w}}
\nc{\fz}{{\mathfrak z}}
\nc{\fx}{{\mathfrak x}}
\nc{\fy}{{\mathfrak y}}

\nc{\al}{{\alpha }}
\nc{\be}{{\beta }}
\nc{\ga}{{\gamma }}
\nc{\de}{{\delta }}
\nc{\del}{{\partial }}
\nc{\ep}{{\varepsilon }}
\nc{\vap}{{\epsilon }}

\nc{\ze}{{\zeta }}
\nc{\et}{{\eta }}
\rc{\th}{{\theta }}

\nc{\vth}{{\vartheta }}

\nc{\io}{{\iota }}
\nc{\ka}{{\kappa }}
\nc{\la}{{\lambda }}
\nc{\vrho}{{\varrho}}
\nc{\si}{{\sigma }}
\nc{\ups}{{\upsilon }}
\nc{\vphi}{{\varphi }}
\nc{\om}{{\omega }}

\nc{\Ga}{{\Gamma }}
\nc{\De}{{\Delta }}
\nc{\nab}{{\nabla}}
\nc{\Th}{{\Theta }}
\nc{\La}{{\Lambda }}
\nc{\Si}{{\Sigma }}
\nc{\Ups}{{\Upsilon }}
\nc{\Om}{{\Omega }}

\nc{\Spec}{\on{Spec}}
\nc{\id}{\on{id}}
\nc{\inv}{ ^{-1}}
\nc{\su}{\subset}
\nc{\ot}{\otimes}
\nc{\un}{\underline}
\nc{\ov}{\overline}
\rc{\k}{\Bbbk}

\nc{\shom}{\cH om}
\nc{\colim}[1]{\underset{#1}{\underrightarrow\lim} \; }
\nc{\Gr}{\cG \tn{r}}
\nc{\one}{\mathds{1}}

\nc{\se}{\section}
\nc{\sse}{\subsection}
\nc{\ssse}{\subsubsection}

\title{Geometric Langlands duality and forms of reductive groups}
\author{Vivek Dhand}

\address{Department of Mathematics \\ Michigan State University}
\email{dhand@math.msu.edu}

\begin{document}

\maketitle
\thispagestyle{empty}

\begin{abstract}
The category of perverse sheaves on the affine Grassmannian of a complex reductive group $G$ gives a canonical geometric construction of the split form of the Langlands dual group $\check G_\bZ$ over the integers \cite{MV}.  Given a field $k$, we give a Tannakian construction of the quasi-split forms of $\check G_k$, as well as a construction of the gerbe associated to an inner form of $\check G_k$.

\end{abstract}

\se{Introduction}
Let $G$ be a complex connected reductive group.  Let $\cO = \bC[[z]]$ and $\cK= \bC((z))$.  Define the following sheaves on $Sch/\bC$:
\[ G_\cO = \un{Hom}_{Sch/\bC}(\Spec(\cO),G) \]
\[ G_\cK = \un{Hom}_{Sch/\bC}(\Spec(\cK),G) \]
Let $\Gr = G_\cK / G_\cO$ be the affine Grassmannian of $G$.  It turns out that $G_\cO$ is representable by an affine group scheme over $\bC$, while $G_\cK$ and $G_\cO$ are only ind-representable.
The Satake category $P_{G_\cO}(\Gr,\bZ)$ is the category of $G_\cO$-equivariant perverse sheaves on $\Gr$ with $\bZ$ coefficients.  It turns out that $P_{G_\cO}(\Gr,\bZ)$ is a symmetric monoidal category (i.e.\;tensor category) under fusion/convolution.  Mirkovi\'c and Vilonen \cite{MV} showed that the Satake category is tensor equivalent to $Rep(\check G_\bZ)$, the category of finite dimensional rational representations of the integral split form of the Langlands dual group $\check G_\bZ$. 

Let $k$ be a field.  In order to construct non-split forms of $\check G_k$, we rely on descent theory for neutral Tannakian categories over $k$ (see \cite{DM}, \cite{Sa}).  Given a finite Galois extension $K/k$ with Galois group $\Ga$, a neutral Tannakian category over $k$ is the same as a neutral Tannakian category over $K$ equipped with a descent datum relative to $K/k$.  We begin with the construction of quasi-split forms, i.e.\;those with a Borel subgroup defined over $k$.

\begin{thm}
Let $K/k$ be a finite Galois extension with Galois group $\Ga$, and let $\rho: \Ga \to Out(G)$ be a group homomorphism.  There is a natural descent datum on $P_{G_\cO}(\Gr,K)$ relative to $K/k$ such that the corresponding $k$-linear category is tensor equivalent to $Rep(\check G^\rho_k)$, where $\check G^\rho_k$ is a quasi-split form of $\check G_k$ corresponding to $\rho$.
\end{thm} 

The classification of inner forms is more difficult, since the problem is just as complicated as the classification of $\check G_k$ torsors.  Given an inner form $\tilde G_k$ of $\check G_k$, we construct a (non-neutral) Tannakian category which is equivalent to $Rep(\tilde G_k)$.  Let $Br(k)$ denote the Brauer group of $k$ and let $\mu: \pi_1(G) \to Br(k)$ be a group homomorphism. Given $\al \in \pi_1(G) = \pi_0(\Gr)$, let $\Gr_\al$ denote the corresponding component of $\Gr$.  Let $A_\al$ be a central simple $k$-algebra such that $[A_\al] = \mu(\al)$.  We define a tensor category $P(\Gr,\mu)$ of ``perverse sheaves with coefficients in $\mu$."  Roughly speaking, over the component $\Gr_\al$ this category consists of right $A_\al$-modules in the category of perverse sheaves on $\Gr_\al$.  The tensor structure on $P(\Gr,\mu)$ is induced from the tensor structure on the Satake category.
 
\begin{thm}  
Let  $\tilde G_k$ be an inner form of $\check G_k$.  There exists a group homomorphism $\mu: \pi_1(G) \to Br(k)$ such that $P(\Gr,\mu)$ and $Rep(\tilde G_k)$ are equivalent as tensor categories. \end{thm}

 The layout of the paper is as follows:  In section \ref{forms}, we recall the classification of forms of $\check G_k$ in terms of the Galois cohomology of $k$. In section \ref{formalism}, we explain the formalism of descent of Tannakian categories. In sections \ref{qs} and \ref{inner} we prove the main theorems about quasi-split forms and inner forms, respectively.

\sse*{Acknowledgments} This paper is adapted from my doctoral thesis, completed at Northwestern University.  I would like to thank my advisor Kari Vilonen for his insight and guidance.  Thanks are also due to David Nadler, Matt Emerton, Robert Kottwitz, and David Treumann for many useful conversations.

\se{Classification of forms} \label{forms}
Let $\Ga$ be a finite group and let $A$ be a $\Ga$-group, i.e.\;a discrete (not necessarily abelian) group on which $\Ga$ acts by group homomorphisms.  We define
\[ H^0(\Ga,A) = A^\Ga \]
where $A^\Ga$ is the subgroup of elements of $A$ fixed by $\Ga$. We define the set of $A$-valued 1-cocycles to be:
\[ Z^1(\Ga,A) = \{ f:\Ga \to A \mid f(st) = f(s) s(f(t)) \textnormal{ for all } s,t\in \Ga \}  .\]
We say that $a,b \in Z^1(G,A)$ are {\em equivalent} if there exists $c \in A$ such that \[b(s) = c^{-1} a(s) s(c)\] for all $s \in \Ga$.  The set of equivalence classes of 1-cocycles is denoted $H^1(\Ga,A)$.  It is not a group in general, but rather a pointed set whose basepoint is the constant cocycle at the identity of $A$. See Section 1.3 of \cite{PR} for more on non-commutative cohomology.

Let $K/k$ be a finite Galois extension with Galois group $\Ga$. Let $X$ and $Y$ be affine $k$-schemes. We say that $Y$ is a $K/ k$-{\em form} of $X$ if there exists an isomorphism of $K$-schemes $X_K \xrightarrow{\simeq} Y_K$, where $X_K = X \times_{\Spec(k)} \Spec(K)$.

Let $Aut(X_K)$ denote the automorphism group of $X_K$ in the category of $K$-schemes.  Then $\Ga$ acts on $Aut(X_K)$ as follows: given $\ga \in \Ga$ and $\phi \in Aut(X_K)$, the morphism $\ga(\phi)$ is given the the following commutative diagram:
\[ \xymatrix{ X_K \ar@{-->}[rr]^{\ga(\phi)} \ar[d]_{\phi} & & X_K \ar[d]^{\id\times\ga} \\ X_K \ar[r]^{\id \times \ga} & X_K & X_K \ar[l]_{\phi} } \]

Let $\cA(X, K/k)$ denote the set of isomorphism classes of $K/ k$-forms of $X$.  It is a pointed set with basepoint given by $X$. 
\begin{prop} Let $X$ be an affine $k$-scheme.  Let $K/k$ be a finite Galois extension with $\Ga = Gal(K/k)$.  Then there is a natural bijection of pointed sets
\[  \cA(X, K /k) \xrightarrow{\simeq} H^1(\Ga, Aut(X_K)) . \]
\end{prop}

\begin{proof}
See \cite{Sp2}, Prop. 11.3.3 and Prop. 12.3.2.
\end{proof}

This result is also true if $X$ is an affine $k$-group.  The proof is the same except we replace $k$-algebras with Hopf algebras over $k$.  For an affine $k$-group $G$, we define
\[ H^1(k, G) = \colim{K/k \tn{ finite Galois}} H^1(Gal(K/k), G(K)) \]
This is nothing but the first \'etale cohomology of $\Spec(k)$ with coefficients in the sheaf of groups $Aut(X)$.  It classifies the $k$-{\em forms} of $X$, i.e. $k$-schemes which are isomorphic to $X$ over a separable closure $k^{sep}$ of $k$. 

Given a short exact sequence of affine $k$-groups:
\[ 1 \to A \to B \to C \to 1\]
there is a long exact sequence in cohomology:
\[ \dots \to H^0(k,C) \to H^1(k,A) \to H^1(k,B) \to H^1(k,C). \] 
If $A$ is commutative we can extend this sequence to include the boundary map:
\[ \dots \to H^1(k,C) \to H^2(k,A). \]

For the rest of this section, $G$ will be a complex connected reductive group and $\check G_k$ will be the connected split reductive $k$-group with root data dual to that of $G$.  We are interested in calculating $H^1(k, Aut(\check G_k))$.  Consider the short exact sequence
\[ 1 \to Int(\check G_k) \to Aut(\check G_k) \to Out(\check G_k) \to 1 \] 
The forms of $\check G_k$ that lie in the image of the map \[H^1(k, Int(\check G_k)) \to H^1(k, Aut(\check G_k))\] are called {\em inner forms} of $\check G_k$.  Any form that is not inner is called an {\em outer form}.  A {\em quasi-split form} is one which has a Borel subgroup defined over $k$. Quasi-split forms of $\check G_k$ are classified by 
\[H^1(k,Out(\check G_k)) = H^1(k,Out(G))\]
Since $\check G_k$ is split, the action of $Gal(k^{sep}/k)$ on $Out(\check G_k)$ is trivial.  Therefore
\[ H^1(k, Out(\check G_k)) = Hom(Gal(k^{sep}/k), Out(G))/\sim \]
where $\sim$ denotes the equivalence relation induced by conjugation in $Out(G)$. 

If $Q$ is a quasi-split form of $\check G_k$, the pre-image of $Q$ under the map
\[ H^1(k, Aut(\check G_k)) \to H^1(k, Out(\check G_k)) \]
is equal to $H^1(k, Int(Q))$.  In particular, every form of $\check G_k$ is an inner form of a quasi-split group. In any case, we are reduced to calculating $H^1(k,Int(Q))$, where $Q$ is a quasi-split form of $\check G_k$. In what follows, we continue to work with a split group for notational convenience, but formulas do exist for quasi-split groups as well. Consider the short exact sequence
\[ 1 \to Z(\check G_k) \to \check G_k \to Int(\check G_k) \to 1, \]
which yields an exact sequence of pointed sets
\begin{equation} \label{boundary} \dots \to H^1(k,\check G_k) \to H^1(k,Int(\check G_k)) \to H^2(k, Z(\check G_k)). \end{equation}
The first term is difficult to calculate in general, but descriptions of $H^1(k,\check G_k)$ in terms of $G$ do exist over the real numbers (see \cite{Bv}) and finite extensions of $\bQ_p$ (see \cite {K1, K2}). On the other hand, the last term of equation \ref{boundary} is actually an abelian group, and has a nice description:

\begin{prop}\label{H2}
Let $k$ be a field and let $G$ be a complex connected reductive group.  Let $\check G_k$ denote the split connected reductive $k$-group with root data dual to that of $G$. Then there is a natural isomorphism:
\[ H^2(k,Z(\check G_k)) \simeq Hom(\pi_1(G), Br(k)) \]
\end{prop}
\begin{proof}
Let $K/k$ be a finite Galois extension with Galois group $\Ga$.  It suffices to prove that
\[ H^2(\Ga, Z(\check G_K)) \simeq Hom(\pi_1(G), Br(K/k)) \]
Recall that
\[ Z(\check G_K) = Hom(\pi_1(G), K^\times)  \]
and that $\Ga$ acts trivially on $\pi_1(G)$ and in the natural way on $K^\times$. We define a map
\[ Z^2(\Ga, Hom(\pi_1(G), K^\times)) \to Hom(\pi_1(G), Z^2(\Ga, K^\times)) \]
\[ \zeta \mapsto \mu \]
by the formula
\[ \zeta(a,b)(\al) = \mu(\al)(a,b) \]
for all $a, b \in \Ga, \al \in \pi_1(G)$.  This map is clearly a bijection.  In fact, the same map restricts to a bijection
\[ B^2(\Ga, Hom(\pi_1(G), K^\times)) \to Hom(\pi_1(G), B^2(\Ga, K^\times)) \]
Suppose we have $\zeta' \in Z^2(\Ga, Hom(\pi_1(G),K^\times))$ and a map $f: \Ga \to Hom(\pi_1(G), K^\times)$ such that
\[ \zeta(a,b)(\al) - \zeta'(a,b)(\al) = f(a)(\al) + a(f(b))(\al) - f(ab)(\al) \]
Define $g(\al): \Ga \to K^\times$ by the formula
\[ g(\al)(a) = f(a)(\al) \]
Since $\Ga$ acts trivially on $\pi_1(G)$, we see that
\[ a(f(b))(\al) = a(f(b)(\al)) = a(g(\al)(b)) = a(g(\al))(b)  \]
Therefore
\[ \mu(\al)(a,b) - \mu'(\al)(a,b) = g(\al)(a) + a(g(\al))(b) - g(\al)(ab) \]
This implies that $\mu(\al) - \mu'(\al) \in B^2(\Ga, K^\times)$ and so we get an isomorphism on the level of cohomology groups.
\end{proof}

Since $\pi_1(G)$ can be read off from the root datum of $\check G_k$, we are left with the problem of calculating $Br(k)$. For one-dimensional fields, we can use the main theorems of local and global class field theory.  If $k$ is a non-archimedean local field, we have a natural isomorphism
\[ Br(k) \cong \bQ/ \bZ \]
given by taking the Hasse invariant of a central simple algebra.  If $k$ is a global field, we have a short exact sequence
\[ 0 \to Br(k) \to \bigoplus_{v} Br(k_v) \to \bQ / \bZ \to 0\]
where $v$ is a valuation, $k_v$ is the completion with respect to $v$, and the map to $\bQ/\bZ$ is the sum of the local Hasse invariants.  Along with the Hasse principle, these results show that the problem of calculating the forms of reductive groups over a number field can be reduced to class field theory.

\se{Descent formalism}\label{formalism}
In this section $K/k$ will be a finite Galois extension with Galois group $\Ga$.  Let $\cC$ be a $K$-linear category.  A {\em descent datum} on $\cC$ relative to $K/k$ is the following (see \cite{DM}, 3.13):
\begin{enumerate}
\item for each $a \in \Ga$, an equivalence $\be_a : \cC \xrightarrow{\simeq} \cC$,
\item for each $a, b \in \Ga$, a natural isomorphism $\mu(a,b) : \be_{ab} \xrightarrow{\simeq} \be_a \be_b$,
\end{enumerate}
such that 
\begin{enumerate}
\item for every $\la \in K$, $f \in Mor(\cC)$, we have $\be_a(\la f) = a(\la) \be_a(f)$,
\item for every $a,b,c \in \Ga$, $x \in \cC$, the following diagram commutes:
\begin{equation} \label{descent} \xymatrix{  \be_{a b c}(x) \ar[rr]^{\mu(a,bc)_x} \ar[d]_{\mu(ab,c)_x} && \be_a\be_{bc}(x) \ar[d]^{\be_a( \mu(b,c)_x)} \\ \be_{ab}\be_c(x) \ar[rr]_{\mu(a,b)_{\be_c(x)}} && \be_a\be_b\be_c(x)}  \end{equation} 
\end{enumerate}
By abuse of notation, we will write $(\be_a,\mu(a,b))$ for any descent datum on $\cC$.

A $k$-linear category $\cC$ gives rise to a $K$-linear category $\cC_K$ with descent datum.  The objects of $\cC_K$ are $K$-modules in $\cC$, i.e. pairs $(x,\rho)$ with $x \in \cC$ and $\rho: K \to End_\cC(x)$ a morphism of $k$-algebras.  The descent datum on $\cC_K$ is given by
\begin{equation}\label{split_action}  \be_a(x,\rho) = (x, \rho \circ a\inv),  \;\;\; \mu(a,b)_x = \id_x. \end{equation} 

Conversely, a $K$-linear category $\cC$ with descent datum $(\be_a,\mu(a,b))$ gives rise to a $k$-linear category $\cC^\Ga$ such that $(\cC^\Ga)_K \simeq \cC$. 
The objects of $\cC^\Ga$ are of the form $(x,\{a_x\}_{a \in \Ga})$ where $x \in \cC$ and $a_x : x \xrightarrow{\simeq} \be_a(x)$ is an isomorphism such that the following diagram commutes for all $a,b \in \Ga$:
\begin{equation} 
\xymatrix{  x \ar[rr]^{a_x}  \ar[d]_{(ab)_x} & & \be_a(x) \ar[d]^{\be_a(b_x)} \\ \be_{ab}(x) \ar[rr]_{\mu(a,b)_x} & & \be_a\be_b(x) } 
\end{equation}
For ease of notation, we will often write $(x,a_x)$ instead of $(x, \{a_x\}_{a \in \Ga})$. A  morphism in $\cC^\Ga$ from $(x,a_x)$ to $(y,a_y)$ is a morphism $f:x \to y$ in $\cC$ such that the following diagram commutes for all $a \in \Ga$:
\begin{equation} 
\xymatrix{x \ar[rr]^{a_x} \ar[d]_{f} && \be_a(x) \ar[d]^{\be_a(f)} \\ y \ar[rr]_{a_y} && \be_a(y) } 
\end{equation}  
If we wish to make the descent datum explicit, we will write $\cC^{(\be_a,\mu(a,b))}$ instead of $\cC^\Ga$.  

\begin{example}\label{vect}
Let $k$ be a field and let $Vect_k$ denote the category of finite dimensional vector spaces over $k$.  If $K$ is a finite extension of $k$, then $(Vect_k)_K = Vect_K$.  If $K/k$ is a finite Galois extension with Galois group $\Ga$, then the descent datum on $Vect_K$ is given by equation \ref{split_action} and denoted $(\ga_a, \id)$.  Moreover, we have an exact equivalence of $k$-linear abelian tensor categories $Vect_K^{(\ga_a, \id)} \simeq Vect_k$.
\end{example}

\begin{rmk} \label{space_example}
Let $X$ be a topological space and let $a, b, c \in Aut(X)$.  Let $\xi(a,b) : (ab)^* \simeq b^* a^*$ denote the natural isomorphism.  Then the following diagram commutes for any sheaf $\cF$ on $X$:
\[ \xymatrix{  (abc)^* \cF \ar[rr]^{\xi(ab,c)_\cF} \ar[d]_{\xi(a,bc)_\cF} && c^* (ab)^* \cF \ar[d]^{c^*\xi(a,b)_\cF} \\ (bc)^* a^* \cF \ar[rr]^{\xi(b,c)_{a^*\cF}} && c^* b^* a^* \cF }  \]
In particular, suppose we have a finite group $\Ga$ acting on $X$ via a group homomorphism $\rho: \Ga \to Aut(X)$.  Let $\be_a = \rho(a\inv)^*$ and define $\mu(a,b) = \xi(\rho(b\inv), \rho(a\inv)) : \be_{ab} \simeq \be_a \be_b$.  Then the isomorphisms $\mu(a,b)$ satisfy diagram \ref{descent} for any sheaf $\cF$ on $X$.
\end{rmk}

\begin{rmk}\label{sheaf}
Let $X$ be topological space and let $Sh(X,K)$ denote the category of $K$-sheaves on $X$.  This category has a natural descent datum relative to $K/k$.  For any sheaf $\cF \in Sh(X,K)$ and any open set $U \su X$, we define
\begin{equation} (\tilde\ga_a\cF)(U) = \ga_a(\cF(U)) \end{equation}
where $\ga_a$ is the functor defined in example \ref{vect}.  Since $\ga_a$ is exact, we see that $\tilde\ga_a(\cF)$ is a sheaf.  In fact, the functor $\tilde \ga_a $ is also exact because that the natural map of $K$-vector spaces
\[ (\tilde \ga_a \cF)_x = \colim{x \in U} \ga_a(\cF(U)) \to \ga_a( \cF_x  )  \]
is a bijection on the underlying $k$-vector spaces for any $x \in X$. We claim that $\tilde\ga_a(\la f) = a(\la) f$ for any $\la \in K$ and any morphism $f$ in $Sh(X,K)$.  Indeed, we have
\[ \tilde\ga_a(\la f)_U = \ga_a(\la f_U) = a(\la) f_U. \]
We also make the identification $\tilde \ga_{ab} = \tilde \ga_a \tilde \ga_b$.  Therefore $(\tilde \ga_a, \id)$ is a descent datum on $Sh(X,K)$ such that $Sh(X,K)^\Ga = Sh(X,k)$.

Now suppose we are given a continuous action $\rho: \Ga \to Aut(X)$ of $\Ga$ on $X$.  The functors $\rho(a\inv)^*$ and $\tilde \ga_b$ commute for all $a,b\in\Ga$.  Indeed, since $\rho(a\inv)$ is a homeomorphism, we have
\[ (\rho(a\inv)^*\cF)(U) = \cF(a\inv U) \]
and therefore,
\begin{eqnarray*} (\rho(a\inv)^* \tilde \ga_b \cF)(U) &= &  (\tilde \ga_b\cF)(a\inv U)) \\  &=& \ga_b (\cF(a \inv U)) \\ &=& \ga_b((\rho(a\inv)^* \cF)(U) ) \\ &=& (\tilde\ga_b \rho(a\inv)^*\cF)(U).
\end{eqnarray*}
These facts allow us to define a new descent datum $(\be_a, \mu(a,b))$ on $Sh(X,K)$, where
\begin{equation} \be_a = \rho(a\inv)^* \tilde \ga_a = \tilde \ga_a \rho(a\inv)^*,\;\;\;\; \mu(a,b) = \xi(\rho(b\inv), \rho(a\inv)) \end{equation}
\end{rmk}

\begin{rmk}\label{perverse}
Let $(X,\cS)$ be a Whitney stratified complex algebraic variety with an action of $\Ga = Gal(K/k)$.  The functors $\tilde \ga_a$ extend to the bounded derived category $D^b_\cS(X,K)$ of $S$-constructible $K$-sheaves on $X$.  In fact, $\tilde \ga_a$ is $t$-exact with respect to the perverse $t$-structure $(D^{\leq 0}, D^{\geq 0})$ on $D^b_\cS(X,K)$.  To see this, recall that
\[ \cF \in D^{\leq 0} \iff \dim supp(H^j(\cF)) \leq -j \tn{ for all } j \in \bZ\]
\[ \cF \in D^{\geq 0} \iff \dim supp(H^j(\bD \cF)) \leq -j \tn{ for all } j \in \bZ\]
Since $\tilde \ga_a$ is exact, we see that $H^j(\tilde \ga_a \cF) = \tilde \ga_a H^j(\cF)$, and $supp(H^j(\tilde \ga_a \cF)) = supp(\tilde \ga_a H^j(\cF))$.  Since $\tilde \ga_a$ is a tensor functor, we know that it commutes with internal hom, hence $\tilde \ga_a \bD \cF \simeq \bD \tilde \ga_a \cF$.  This shows that
\[ \tilde \ga_a (D^{\leq 0}) \su D^{\leq 0} \tn{ and } \tilde \ga_a (D^{\geq 0}) \su D^{\geq 0} \]
Therefore, we get exact functors, for all $a \in \Ga$:
\[ {}^p \tilde \ga_a : P_\cS(X,K) \simeq P_{\cS}(X,K) \]
where $P_{\cS}(X,K) = D^{\leq 0} \cap D^{\geq 0}$ denotes the abelian category of $\cS$-constructible perverse $K$-sheaves on $X$.  If no confusion will result, we will continue to denote ${}^p \tilde \ga_a$ as simply $\tilde \ga_a$.

If $\Ga$ acts on $X$ via $\rho: \Ga \to Aut(X)$, then we have exact functors, for all $b \in \Ga$:
\[ \rho(b \inv)^* : P_\cS(X,K) \simeq P_\cS(X,K) \]
which commute with all the functors $\tilde \ga_a$.  Therefore, 
\[(\rho(b \inv)^* \tilde \ga_a, \;\; \xi(\rho(b\inv), \rho(a\inv)))\]
is a descent datum on $P_\cS(X,K)$.
\end{rmk}

Now suppose we have $K$-linear categories $\cC$ and $\cC'$ with descent data $(\be_a,\mu(a,b))$ and $(\be'_a,\mu'(a,b))$, respectively. A functor $F: \cC \to \cC'$ is {\em compatible} with the descent data if there exist natural isomorphisms $ \eta(a): F \be_a \xrightarrow{\simeq} \be'_a F $, for each $a \in G$, such that the following diagram commutes:
\begin{equation} \label{eta}
\xymatrix{F \be_{ab} x \ar[rr]^{F(\mu(a,b)_x)} \ar[d]_{\eta(ab)_x }  & & F \be_a \be_b x \ar[rr]^{\eta(a)_{\be_b x}}   & & \be_a' F \be_b x  \ar[d]^{\be_a' (\eta(b)_X)}  \\  \be'_{ab} Fx \ar[rrrr]^{\mu'(a,b)_{Fx}} &&  && \be_a'\be_b' Fx} 
\end{equation}
In addition, if $F$ is an equivalence we say that $(\be_a,\mu(a,b))$ and $(\be'_a,\mu'(a,b))$ are equivalent.

\begin{prop} A functor $(\cC \xrightarrow{F} \cC')$ that is compatible with descent data induces a functor $F^\Ga: \cC^\Ga \to \cC'^\Ga$. \end{prop}
\begin{proof}
We define $F^\Ga$ on objects of $\cC^\Ga$ as follows:
\[ x \mapsto Fx \]
\[ a_x \mapsto \eta(a)_x F(a_x) \]
We claim that this defines an object of $\cC'^\Ga $.  Consider the following diagram:
\[ \xymatrix{  Fx \ar[rr]^{F(a_x)} \ar[d]_{F((ab)_x)} & & F\be_a x \ar[rr]^{\eta(a)_x} \ar[d]^{F \be_a (b_x)} & & \be_a' Fx \ar[d]^{\be_a' F (b_x)}   \\  F \be_{ab} x \ar[rr]_{F(\mu(a,b)_x)} \ar[d]_{\eta(ab)_x }  & & F \be_a \be_b x \ar[rr]_{\eta(a)_{\be_b x}}   & & \be_a' F \be_b x  \ar[d]^{\be_a' (\eta(b)_x)}  \\  \be'_{ab} Fx \ar[rrrr]^{\mu'(a,b)_{Fx}} &&  && \be_a'\be_b' Fx  } \]
The square on the top left commutes because $(x,a_x)$ is an object of $\cC^\Ga$.  The top right square commutes because $\eta(a)$ is a natural transformation. The bottom pentagon commutes by equation \ref{eta}.  Therefore, the compositions along the outside edges are equal:
\[  \be_a' ( \eta(b)_x F(b_x) )  \; \eta(a)_x F(a_x) = \mu'(a,b)_{Fx}  \; \eta(ab)_x F((ab)_x) \]
which establishes the claim.  

We define $F^\Ga(f) = F(f)$ on morphisms and check that the following diagram commutes:
\[  \xymatrix{Fx \ar[rr] ^{F(a_x)}\ar[d]_{F(f)} & & F \be_a x \ar[rr]^{\eta(a)_x}  \ar[d]^{F \be_a (f)} & & \be_a' F x \ar[d]^{\be_a' F (f)}  \\ F y \ar[rr]_{F(a_y)} & & F \be_a y \ar[rr]_{\eta(a)_y} & & \be_a' F y  }\]
The left square commutes because $f$ is a morphism in $\cC^\Ga$ and the right square commutes because $\eta(a)$ is a natural transformation.  Therefore $F^\Ga(f)$ is a morphism in $\cC'^\Ga$. The rest of the proposition follows easily.
\end{proof}

We collect some basic facts about descent of categories:
\begin{prop} \label{descent_facts}
Let $K/k$ be a finite Galois extension with Galois group $\Ga$.  Let $\cC$ be a $K$-linear category with descent datum $(\be_a, \mu(a,b))$ relative to $K/k$.
\begin{enumerate}
\item If $\cC$ is abelian and $\be_a$ is exact and additive for every $a \in \Ga$, then $\cC^\Ga$ is abelian and the forgetful functor $\cC^\Ga \to \cC$ is exact.
\item If $\cC$ is a tensor category, $\be_a$ is a tensor functor for every $a \in \Ga$, and $\mu(a,b)$ is a morphism of tensor functors for every $a,b \in \Ga$, then $\cC^\Ga$ is a tensor category.
\item Suppose that $\cC$ is an abelian tensor category equipped with a tensor descent datum, i.e. a descent datum satisfying the conditions in (1) and (2) above.  If $\om: \cC \to Vect_K$ is an exact faithful tensor functor which is compatible with the tensor descent datum, then $\om^\Ga : \cC^\Ga \to Vect_k$ is an exact faithful tensor functor.
\end{enumerate}
\end{prop}

\begin{proof}
(1) Since each $\be_a$ is additive, we have $0 \simeq \be_a 0$, and this defines a zero object in $\cC^\Ga$.  The direct sum of $(x, \{a_x\})$ and $(y, \{a_y\})$ in $\cC^\Ga$ is defined to be:
\[ (x \oplus y, \{a_x \oplus a_y\}) \]
This defines an object of $\cC^\Ga$ since $\be_a$ is additive for all $a \in \Ga$, and it clearly satisfies the desired universal property. Let $f : (x, \{a_x\}) \to (y, \{a_y\}) $ be a morphism in $\cC^\Ga$.  Since $\be_a$ is exact, we obtain the following commutative diagram with exact rows:
\[ \xymatrix{ 0 \ar[r] & Ker(f) \ar@{..>}[d]^{\ov{a_x}} \ar[r] &  x \ar[d]^{a_x} \ar[r] & y \ar[d]^{a_y}  \\ 0 \ar[r] & Ker( \be_a f) \ar@{..>}[d]^{i_a(f)} \ar[r] & \be_a x \ar[r] \ar@{=}[d] & \be_a y \ar@{=}[d]  \\ 0 \ar[r] & \be_a(Ker(f)) \ar[r] & \be_a x \ar[r] & \be_a y  } \]
The composition of the vertical dotted arrows defines $a_{Ker(f)} : Ker(f) \simeq \be_a(Ker(f))$.  The following commutative diagram implies that kernels exist in $\cC^\Ga$:
\[ \xymatrix{ Ker(f) \ar[r]^{\ov{a_x}} \ar[d]_{\ov{(ab)_x}} & Ker(\be_a f) \ar[r]^{i_{a}(f)} \ar[d]^{\ov{\be_a(b_x)}} & \be_a(Ker(f)) \ar[d]^{\be_a(\ov{b_x})}  \\ Ker(\be_{ab} f) \ar[r]^{\ov{\mu(a,b)_x}} \ar[d]_{i_{ab}(f)}  & Ker(\be_a \be_b f) \ar[r]^{i_a(\be_b f)} & \be_a(Ker(\be_b f))\ar[d]^{\be_a(i_b(f))}  \\ \be_{ab} Ker(f) \ar[rr]^{\mu(a,b)_{Ker(f)}}  & & \be_a \be_b Ker(f) \\  } \]
The bottom pentagon commutes because of the identities between various short exact sequences coming ``out of the page." A similar argument shows that cokernels exist in $\cC^\Ga$.  The exactness axiom and the exactness of the forgetful functor $\cC^\Ga \to \cC$ follow immediately.
\medskip

(2) Let $(\cC, \otimes,\one,\phi,\psi)$ be a tensor category with unit object $\one$, associativity constraint $\phi$, and commutativity constraint $\psi$. We define the tensor structure on $\cC^\Ga$ as follows:
\[(x,\{a_x\})\otimes (y,\{a_y\})=(x \otimes y, a_{x\otimes y} = \ka(a)_{x,y} \circ (a_x \otimes a_y))\]  
where 
\[\ka(a)_{x,y} : \be_a(x) \otimes \be_a(y) \simeq \be_a(x \otimes y)  \] 
is the natural isomorphism making $\be_a$ into a tensor functor.  

This formula gives an object of $\cC^\Ga$ because the following diagram commutes:
\[ \xymatrix{ x \otimes y \ar[r]^{a_x \otimes a_y} \ar[d]_{(ab)_x \otimes (ab)_y} & \be_a x \otimes \be_a y \ar[r]^{\ka(a)_{x,y}} \ar[d]^{\be_a(b_x) \otimes \be_a(b_y)} & \be_a(x \otimes y) \ar[d]^{\be_a(b_x \otimes b_y)}  \\ \be_{ab} x \otimes \be_{ab} y \ar[r] \ar@{}[rd]^{\mu(a,b)_x \otimes \mu(a,b)_y} \ar[d]_{\ka(ab)_{x,y}}  & \be_a \be_b x \otimes \be_a \be_b y \ar[r] \ar@{}[rd]^{\ka(a)_{\be_b x, \be_b y}} & \be_a(\be_b x \otimes \be_b y) \ar[d]^{\be_a(\ka(b)_{x,y})}  \\ \be_{ab}( x\otimes y) \ar[rr]_{\mu(a,b)_{x\otimes y}}  & & \be_a \be_b (x \otimes y) \\  } \]
The bottom pentagon commutes because $\mu(a,b)$ is a morphism of tensor functors, and because $\be_a(\ka(b)_{x,y})\circ \ka(a)_{\be_b x, \be_b y}$ is the compatibility isomorphism of the composite tensor functor $\be_a \be_b$.

Note that there is a unique isomorphism $\one \simeq \be_a \one$ since $\be_a$ is a tensor functor.  Therefore we get a well defined unit object in $\cC^\Ga$.  Finally, we need to check that $\phi$ and $\psi$ are morphisms in $\cC^\Ga$.  The case of the commutativity constraint follows from this commutative diagram:
\[ \xymatrix{ x\otimes y \ar[rr]^{a_x\otimes a_y} \ar[d]^{\psi_{x,y}} && \be_a x \otimes \be_a y \ar[rr]^{\ka(a)_{x,y}} \ar[d]^{\psi_{\be_a x, \be_a y}} &&  \be_a(x \otimes y)\ar[d]^{\be_a(\psi_{x,y})} \\ y \otimes x \ar[rr]^{a_y\otimes a_x} && \be_a y \otimes \be_a x \ar[rr]^{\ka(a)_{y,x}} &&  \be(y \otimes x)}  \]
Similarly we can set up the following commutative diagram for the associativity constraint:
\[ \xymatrix{ x \!\otimes \!(y \!\otimes\! z) \ar[rr]^<<<<<<<<<{a_x \otimes (a_y\otimes a_z)} \ar[d]_{\phi_{x,y,z}}  && \be_a x \!\otimes\! (\be_a y \!\otimes\! \be_a z) \ar[r] \ar@{}[rd]^<<<<<<<<{1\otimes \ka(a)_{y,z}} \ar[d]_{\phi_{\be_a x,\be_a y, \be_a z}} & \be_a x \!\otimes\! \be_a (y \!\otimes\! z) \ar[r] \ar@{}[rd]^<<<<<<<<{\ka(a)_{x,y\otimes z}} & \be_a (x \!\otimes\! (y \!\otimes\! z)) \ar[d]^{\be_a(\phi_{x,y, z})} \\ (x \!\otimes\! y)\! \otimes\! z \ar[rr]_<<<<<<<<<{(a_x \otimes a_y) \otimes a_z}  && (\be_a x \!\otimes\! \be_a y) \!\otimes\! \be_a z \ar[r]\ar@{}[ru]_<<<<<<<<{\ka(a)_{x,y}\otimes 1} & \be_a (x \!\otimes\! y) \!\otimes\! \be_a z \ar[r]\ar@{}[ur]_<<<<<<<<{\ka(a)_{x\otimes y,z}} & \be_a((x \!\otimes\! y) \!\otimes\! z )  } \]

(3)  We have already shown that $\om^\Ga : \cC^\Ga \to Vect_k$ is a functor. The proof that $\om^\Ga$ is a tensor functor is tedious but straightforward.  The isomorphism $\om^\Ga((x,a_x)\otimes(y,a_y)) \simeq \om^\Ga(x,a_x) \otimes \om^\Ga(y,a_y)$ comes from the isomorphism $\om(x \otimes y) \simeq \om(x) \otimes \om(y)$. To show that $\om$ is exact, let
\[ 0\to (x,a_x) \to (y,a_y) \to (z,a_z) \to 0 \]
be a short exact sequence in $\cC^\Ga$. Since $\cC^\Ga \to \cC$ is exact, we see that
\[ 0 \to x \to y \to z \to 0 \]
is exact in $\cC$.  Since $\om$ is exact, we get that
\[ 0 \to \om(x) \to \om(y) \to \om(z) \to 0 \]
is exact.
Since the forgetful functor $Vect_K^\Ga \to Vect_K$ is exact, we conclude that
\[0 \to (\om(x), \eta(a)_x \om(a_x)) \to(\om(y), \eta(a)_y \om(a_y)) \to (\om(z), \eta(a)_z \om(a_z)) \to 0\]
is exact.  Finally, if $ \om^\Ga(x, a_x) = (\om(x), \eta(x) \om(a_x)) = 0$, then $\om(x) = 0$, which implies that $x=0$ because $\om$ is faithful.

\end{proof}

\begin{cor} \label{Tannakian}
Let $K/k$ be a finite Galois extension with Galois group $\Ga$. Let $\cC$ be a neutral Tannakian category over $K$ with a fibre functor $\om : \cC \to Vect_K$.  Suppose we have a tensor descent datum on $\cC$ relative to $K/k$ (see Proposition \ref{descent_facts}).  Then $\cC^\Ga$ is a neutral Tannakian category over $k$ with fibre functor $\om^\Ga : \cC^\Ga \to Vect_k$.
\end{cor}
\begin{proof}
We have shown every part of this corollary except that $\cC^\Ga$ is rigid.  Let $Dx$ denote the dual of an object $x \in \cC$ and let $i_x : x \simeq DDx$ be the duality isomorphism.  We also have isomorphisms $\phi(a)_x : D(\be_a x) \simeq \be_a(Dx)$ (\cite{DM}, Proposition 1.9).  We define the dual of $(x,a_x)$ to be $(Dx, \phi(a)_x \circ (Da_x)\inv)$.  The commutativity of the following diagram shows that this is an object of $\cC^\Ga$:
\[\xymatrix{ Dx & D(\be_a x) \ar[l]_{Da_x}  \ar[r]^{\phi(a)_x}  & \be_a(Dx) \\ D(\be_{ab} x) \ar[u]^{D((ab)_x)} \ar[d]_{\phi(ab)_x} & D(\be_a\be_b x) \ar[l]^{D(\mu(a,b)_x)} \ar[r]_{\phi(a)_{\be_b x}} \ar[u]_{D(\be_a(b_x))}  & \be_a( D \be_b x) \ar[u]_{\be_a(Db_x)} \ar[d]^{\be_a(\phi(b)_x)} \\ \be_{ab}(Dx) \ar[rr]^{\mu(a,b)_{Dx}} & & \be_a\be_b (Dx) }\]
It remains to check that $(Dx, \phi(a)_x \circ (Da_x)\inv)$ is the dual of $(x, a_x).$  This follows from the following commutative diagram, for any $f : x \to \one$:
\[ \xymatrix{  Dx \ar[d]_{Df} & D \be_a x \ar[l] \ar[d]^{D(\be_a f)} \ar[r]^{\phi(a)_x} & \be_a Dx \ar[d]^{\be_a (Df)} \\ D\one & D\be_a \one \ar[l]_\simeq \ar[r]^\simeq & \be_a D \one } \]
where we identify $\one = D\one$ and $\one \simeq \be_a \one$ is the unique isomorphism.
\end{proof}

\se{Quasi-split groups}\label{qs}

Let $G$ be a complex connected reductive algebraic group. Let $\cO = \bC[[z]]$ be the ring of formal power series in one variable, and let $\cK = Frac(\cO) = \bC((z))$ be its fraction field.   We think of $G(\cO)$ as the complex points of a group scheme $G_\cO$.  Similarly, $G(\cK)$ has the structure of an ind-group scheme $G_\cK$, and $G(\cK)/G(\cO)$ has the structure of an ind-scheme $\Gr$ called the affine Grassmannian of $G$.  Since we work with sheaves in the classical topology, we can think of $G_\cO, G_\cK,$ and $\Gr$ as reduced ind-schemes.

Let $T \su G$ be a maximal torus in $G$ and fix a Borel subgroup $B \supset T$.  Let $N$ denote the unipotent radical of $B$. Let $X^*(T) = Hom(T,\bC^\times)$ and $X_*(T) = Hom(\bC^\times,T)$ denote the weight and coweight lattices of $T$, respectively.  Let $\De \su R \su X_*(T)$ be the set of (positive) roots and $\check \De \su \check R \su X_*(T)$ to set of (positive) coroots with respect to $B$.  Let $\La^+ \su X^*(T)$ (resp. $\La_+ \su X_*(T)$) be the set of dominant weights (resp. coweights), and let $\Pi \su \De$ denote the set of simple roots.  Let $W = N_G(T)/T$ be the Weyl group of $G$.  Let $\fg = Lie(G)$ and let
\[\fg = \fg_0 \oplus(\bigoplus_{\al \in R} \fg_\al)\]
be the adjoint decomposition.  Fix a non-zero vector $x_\al \in \fg_\al$ for each $\al \in \Pi$.  We have a short exact sequence
\[ 1 \to Int(G) \to Aut(G) \to Out(G) \to 1 \]
The surjection $Aut(G) \to Out(G)$ restricts to an isomorphism $Aut(G,B,T,\{x_\al\}_{\al \in \Pi}) \simeq Out(G)$ (\cite{Sp1}, 2.14).  Via this splitting, we conclude that $Out(G)$ acts on each of the following objects $G$, $B$, $T$, $X^*(T)$, $X_*(T)$, $R$, $\check R$, $\De$, $\check \De$, $\La^+$, $\La_+$, and $\Pi$.  Since $T \simeq B/N$, we get an induced action on $N$. In addition, $Out(G)$ also acts on $N_G(T)$.  To see this let $\de \in Aut(G,B,T, \{x_\al\})$ and $g \in N_G(T)$.  For all $t \in T$, we have
\[ Int(\de(g))(t) = \de(g) t \de(g)\inv = \de( g \de\inv(t) g\inv ) = \de( Int(g)(\de\inv(t)) ) \]
Therefore, $\de(g) \in N_G(T)$.  This immediately implies that $Out(G)$ acts on $W$.

Given $\la\in X_*(T) $, we define $\tilde \la \in G(\cK)$ to be the composition of the following maps:
\[ \Spec(\cK) \to \Spec(\bC[z^{\pm 1}]) \xrightarrow{\la} T \su G \]
Let $L_\la = \tilde \la \cdot G_\cO$ denote the image of $\tilde \la$ in $\Gr$.  The left action of $G_\cO$ on $\Gr$ gives us an orbit stratification, denoted by the symbol $\cS$.  We write
\[ \Gr = \bigsqcup_{\la} \Gr^\la \]
where $\la \in \La_+$ and $\Gr^\la = G_\cO \cdot L_\la$.  Moreover, the closure of an orbit looks like
\[ \ov{\Gr^\la} = \bigsqcup_{\mu} \Gr^\mu \]
where $\la-\mu $ is a sum of positive coroots.

For any $\nu \in X$, we have the semi-infinite orbit 
\[S_\nu = N_\cK \cdot L_\nu. \]
These orbits also satisfy closure relations:
\[ \ov{S_\nu} = \bigsqcup_{\eta \leq \nu} S_\eta  \]
where $\nu - \eta$ is a sum of positive coroots.

Let $\de \in Out(G) \simeq Aut(G,B,T, \{x_\al\})$.  Note that $\de$ acts on $G_\cO \su G_\cK$, hence induces an automorphism of $\Gr$, which we also call $\de$.  The action of $\de$ is compatible with the above stratifications:
\[ \de(L_\nu) = L_{\de(\nu)} \]
\[ \de(\Gr^\la) = \Gr^{\de(\la)} \]
\[ \de(S_\nu) = S_{\de(\nu)} \]
for all $\la \in \La_+$, $\nu \in X_*(T)$.

Let $k$ be a field and let $D^b_{\cS}(\Gr,k)$ denote the bounded derived category of $\cS$ constructible $k$-sheaves on $\Gr$.  For example, we could define this category as the colimit of the categories $D^b_{\cS}(\ov{\Gr^\la},k)$ for $\la \in \La_+$.  Let $P_{\cS}(\Gr,k) \su D^b_{\cS}(\Gr,k)$ denote the abelian full subcategory of $\cS$-constructible perverse $k$-sheaves.  We can characterize these sheaves as follows (\cite{MV}, Lemma 3.9): $\cA \in D^b_{\cS}(\Gr,k)$ is perverse if and only if, for all $\nu \in X_*(T)$, we have \[H^k_c(S_\nu,\cA) = 0 \tn{ for } k \neq 2\rho(\nu) \] where $\rho \in X^*(T)$ is equal to half the sum of positive roots.  

Let $P_{G_\cO}(\Gr,k)$ denote the category of $G_\cO$-equivariant perverse sheaves on $\Gr$.  It turns out that $P_{G_\cO}(\Gr,k)$ and $P_\cS(\Gr,k)$ are naturally equivalent (\cite{MV}, Proposition 2.1).  We recall the fusion/convolution tensor structure on $P_{G_\cO}(\Gr,k)$.  Consider the following diagram:
\begin{equation} \label{conv} \xymatrix{ \Gr \times \Gr & G_\cK \times \Gr \ar[r]^<<<<q \ar[l]_{p} & G_\cK \times_{G_\cO} \Gr \ar[r]^<<<<m & \Gr } \end{equation}
Since $G_\cO$ acts freely on $G_\cK \times \Gr$, we have an equivalence \cite{BL}:
\[ q^* : P(G_\cK \times_{G_\cO} \Gr ) \simeq P_{G_\cO}(G_\cK \times \Gr) \]
Let $q_*^{G_\cO}$ denote a quasi-inverse to $q^*$. The convolution of two sheaves $\cA, \cB \in P_{G_\cO}(\Gr, k)$ is defined as follows:
\[ \cA \ast \cB = Rm_* (q_*^{G_\cO} p^*(\cA \boxtimes \cB) ) \]
It is a non-trivial result that the convolution is a perverse sheaf.  It follows from the fact that $m$ is a semi-small map.  The associativity constraint on this monoidal product is essentially formal.  However, the commutativity constraint is more difficult to construct, and involves the fusion product, which must be defined globally on a smooth complex curve $X$.  For our purposes, it suffices to consider $X = \bA^1_\bC$.  Let $\Gr_{X^n} = \Gr \times_{\Spec(\bC)} X^n$.  There is a version of diagram \ref{conv} of ind-schemes over $X$, which defines the fusion product.  In general, a braided product which is compatible with a monoidal product induces a symmetric monoidal structure.  Since fusion is braided and compatible with convolution, we get a tensor structure on $P_{G_\cO}(\Gr,k)$. The global cohomology functor $\bH^*$ is a fibre functor for $P_{G_\cO}(\Gr,k)$.  We recall the main result of \cite{MV}: 

\begin{thm} There is an equivalence of tensor categories:
 \[P_{G_\cO}(\Gr,k) \simeq Rep(\check G_k)\] 
 where $\check G_k$ is the split reductive group scheme over $k$ with root data dual to that of $G$.
\end{thm} 
In the statement of this theorem, the coefficients are allowed to be any commutative ring.  However, for the purposes of this paper we will only consider field coefficients. 

Let us turn now to the case of quasi-split groups. As noted in section \ref{forms}, the quasi-split forms of $\check G_k$ are classified by
\[ H^1(k, Out(\check G_k)) = H^1(k, Out(G)) = Hom(Gal(k^{sep}/k), Out(G))/\sim. \]
Therefore, the data of a quasi-split form of $\check G_k$ can be represented by a group homomorphism from a finite Galois group to $Out(G)$.

\begin{thm}
Let $K/k$ be a finite Galois extension with Galois group $\Ga$, and let $\rho: \Ga \to Out(G)$ be a group homomorphism.  There is a natural descent datum on $P_{G_\cO}(\Gr,K)$ such that the corresponding $k$-linear category $P_{G_\cO}(\Gr,K)^\Ga$ is tensor equivalent to $Rep(\check G^\rho_k)$, where $\check G^\rho_k$ is a quasi-split form of $\check G_k$ corresponding to $\rho$.
\end{thm}

\begin{proof}
Let $\de \in Out(G) \simeq Aut(G,B,T, \{x_\al\})$. We claim that 
\[\de^*: P_{G_\cO}(\Gr,k) \to P_{G_\cO}(\Gr,k)\] 
is a tensor functor.  Let $\cA \in P_{G_\cO}(\Gr,k)$ be a sheaf supported on some $\ov{\Gr^\la}$, such that the action of $G_{\cO}$ on $\ov{\Gr^\la}$ factors through a finite dimensional group $G_{\cO_n}$, where $\cO_n = \cO/(z^n)$.  Since $\cA \in P_\cS(\ov{\Gr^\la},k)$ is $G_{\cO_n}$-equivariant, we have an isomorphism 
\[\phi: p_2^* \cA \simeq a^* \cA \] 
where $a, p_2: G_{\cO_n} \times \ov{\Gr^\la} \to \ov{\Gr^\la}$ are the action map and the projection to the second factor, respectively.  Moreover, the morphism $\phi$ has the property that $i^* \phi = \id$, where $i : \{1\} \times \ov{\Gr^\la} \to \ov{\Gr^\la}$.  We claim that $\de^*\cA$ is also an object of $P_{G_{\cO_n}}(\ov{\Gr^\la},k)$. Indeed, the following commutative diagram:
\[ \xymatrix{ G_{\cO_n} \times \ov{\Gr^\la} \ar[rr]^{p_2, a} \ar[d]_{(\de,\de)} && \ov{\Gr^\la} \ar[d]^\de \\  G_{\cO_n} \times \ov{\Gr^\la} \ar[rr]^{p_2, a}  &&  \ov{\Gr^\la} } \]
induces the following sequence of isomorphisms:
\[ p_2^* \de^* \cA  \simeq  (\de,\de)^*p_2^*\cA  \xrightarrow{(\de,\de)^*\phi} (\de,\de)^*a^* \cA  \simeq a^*\de^* \cA . \]
and clearly, $i^* (\de,\de)^* \phi = \id$.

In order to show that $\de^*$ is compatible with the constraints, we must use the formulation of the tensor structure as a fusion product (\cite{MV}, Section 5).  Let $X = \bA^1_\bC$, and let $\Gr_X$ and $\Gr_{X^2}$ denote the families of Grassmannians over $X$ and $X^2$, respectively.  We can choose a global coordinate on $X$ to trivialize $\Gr_X$ over $X$.  Let $\tau: \Gr_X \to \Gr$ denote the projection.  Let $\De_X \su X^2$ be the diagonal, and let $U = X^2 - \De_X$.  We have that $\Gr_{X^2} |_{\De_X} \cong \Gr_X$ and $\Gr_{X^2} |_U \cong (\Gr_X \times \Gr_X)|_U$.  Let $i: \Gr_X \to \Gr_{X^2}$ and $j:  (\Gr_X \times \Gr_X)|_U \to \Gr_{X^2}$ be the corresponding inclusions.  Let $\tau^o = \tau^*[1]$ and $i^o = i^*[-1]$.  Let
\[ \ast_X : P_{G_{X,\cO}}(\Gr_X,k) \times P_{G_{X,\cO}}(\Gr_X,k)  \to P_{G_{X,\cO}}(\Gr_{X^2},k)   \]
denote the fusion product. Given $\cA, \cB \in P_{G_\cO}(\Gr,k)$, we need to construct an isomorphism
\[  (\de^*\cA) \ast (\de^* \cB) \xrightarrow{c_{\cA, \cB}} \de^*(\cA \ast \cB). \]
Consider the following diagram:
\begin{equation}\label{fusion} \xymatrix{ \tau^o( \de^* \cA \ast \de^* \cB ) \ar[d]_{\simeq} \ar@{-->}[rr]^\simeq  & & \tau^o(\de^*(\cA \ast \cB)) \\ i^o(\tau^o(\de_1^* \cA) \ast_X \tau^o(\de^* \cB)   ) \ar[d]_{\simeq} & & \de^* \tau^o(\cA \ast \cB) \ar[u]_\simeq \\ i^o(j_{!*} {}^pH^0(\tau^o \de^* \cA \overset{L}\boxtimes \tau^o \de^* \cB)|_U) \ar[d]_\simeq && \de^* i^o(\tau^o\cA \ast_X \tau^o \cB)  \ar[u]_\simeq \\  i^o(j_{!*} {}^pH^0((\de,\de)^*(\tau^o \cA \overset{L}\boxtimes \tau^o \cB))|_U) \ar[rr]^\simeq && \de^* i^o(j_{!*} {}^pH^0(\tau^o \cA \overset{L}\boxtimes \tau^o \cB)|_U) \ar[u]_\simeq } \end{equation}
Here we have used equation (5.10) from \cite{MV} along with the fact that $\de$ is an automorphism of $(\Gr_X \times \Gr_X)|_U$, which implies that it commutes with $j_*, j_!$ and $j_{!*}$. Specializing the dotted arrow to any point of $\De_X$ yields the desired isomorphism $c_{\cA,\cB}$.  This argument also shows that the semi-linear functor
\[ \tilde \ga_a : P_{G_\cO}(\Gr,K) \to P_{G_\cO}(\Gr,K) \]
is a tensor functor for all $a \in \Ga$.

Next we define $\be_\de := (\de^{-1})^*: P_{G_\cO}(\Gr,k) \simeq P_{G_\cO}(\Gr,k)$ and note that it is a tensor equivalence.  For $\de_1, \de_2 \in Out(G)$, there is a natural isomorphism $\be_{\de_1 \de_2} \simeq \be_{\de_1} \be_{\de_2}$ described in example \ref{space_example}.  Let us check that the functors $\{\be_\de\}_{\de \in Out(G)}$ are independent of the choice of $\{x_\al \}$.  Any two splittings of
\[ 1 \to Int(G) \to Aut(G) \to Out(G) \to 1 \]
are related by an inner automorphism the form $Int(t)$ for $t \in T(\bC)$ (\cite{Sp1}, 2.14). Therefore, it suffices to show that $\psi_t^* : P_{G_\cO}(\Gr,k) \to P_{G_\cO}(\Gr,k)$ is tensor equivalent to the identity functor for any $t \in T(\bC)$, where $\psi_t : \Gr \to \Gr$ denotes left multiplication by $t$.  Consider the underlying functor $\tilde \psi_t^* : P_\cS(\Gr,k) \to P_\cS(\Gr,k)$.  Any sheaf $\cA \in P_{G_\cO}(\Gr,k)$ is {\em a fortiori} $G$-equivariant, and hence comes equipped with an isomorphism $\cA \simeq \tilde \psi_t^*\cA $ which commutes with morphisms in $P_{G_\cO}(\Gr,k)$.  Using the same argument as in diagram \ref{fusion}, we conclude that $\psi_t^*$ is a tensor functor and that $\id \simeq \psi_t^* $ is an isomorphism of tensor functors.

Let us show that the functors $\{\be_\de\}_{\de \in Out(G)}$ commute with $\bH^*$. Recall that 
\[ \bH^* \simeq \bigoplus_{\nu \in X_*(T)} F_{\nu} \]  
where $F_\nu$ is the weight functor $H^{2\rho(\nu)}_c(S_\nu, \underline{\;\; } )$.  This decomposition yield a canonical split maximal torus $\check T_\k \su \check G_k$:
\begin{equation}\label{Satake} \xymatrix{ P_{G_\cO}(\Gr,k) \ar[rr]^\simeq \ar[d]_{\bigoplus_{\nu} F_\nu} && Rep(\check G_k) \ar[d] \\ Vect_k(X_*(T)) \ar[rr]^\simeq  && Rep(\check T_k) }\end{equation}
where $Vect_k(X_*(T))$ denotes the category of finite dimensional $X_*(T)$-graded vector spaces over $k$.

Consider the following commutative diagram:
\[ \xymatrix{ L_\nu \ar[r]^a \ar[d]^{\de} & S_\nu \ar[r]^b \ar[d]^{\de} & \Gr \ar[d]^{\de} \\ L_{\de(\nu)} \ar[r]^c & S_{\de(\nu)} \ar[r]^d & \Gr }  \]
Since $\de$ is an automorphism, it is proper, so $\de_* = \de_!$.  By uniqueness of adjoints, $\de^* = \de^!$. Therefore, for any sheaf $\cA \in P_{G_\cO}(\Gr,k)$ we have:
\begin{eqnarray*}
 H^{2\rho(\nu)}_c(S_\nu,\de^* \cA) &=& a^! b^* \de^* \cA \\  & = & a^! \de^* d^* \cA \\ & = & a^!\de^!d^* \cA \\ &=&  \de^! c^! d^* \cA \\  &=& H^{2\rho(\de(\nu))}_c (S_{\de(\nu)}, \cA).
\end{eqnarray*}
Therefore, the weight functor $F_\nu$ applied to $\be_\de \cA$ gives:
\[ F_\nu(\be_\de \cA) = F_{\de^{-1}(\nu)}(\cA). \]

Since $\bH^* \simeq \bigoplus_{\nu \in X_*(T)} F_\nu$ after forgetting the coweight grading, we find that
\[  \bH^*(\be_\de \cA) \simeq \bigoplus_{\nu \in X_*(T)} F_{\de^{-1}\nu}(\cA) \simeq \bH^*(\cA). \]

We now change the coefficient ring to $K$. By combining the above discussion with remark \ref{perverse}, we have that 
\[ (\be_{\rho(a)} \tilde \ga_a, \;\; \xi(b\inv,a\inv))\] 
is a tensor descent datum on $P_{G_\cO}(\Gr,K)$ which is compatible with the fibre functor.  By Corollary \ref{Tannakian}, we conclude that $P_{G_\cO}(\Gr,K)^\Ga$ is a neutral Tannakian category.  Let $\check G^\rho_k$ denote the corresponding $k$-form of $\check G_k$.  Since the weight functors are compatible with the descent datum, we get a maximal $k$-torus $\check T^\rho_k \su \check G^\rho_k$. It remains to show that $\check G^\rho_k$ is quasi-split.

Let us recall the ``Pl\"ucker construction"  of the canonical Borel subgroup $\check B_k \supset \check T_k$ (see, for example, the appendix of \cite{Na}): given $\la \in \La_+$, let $j : \Gr^\la \to \Gr $ be the inclusion and let $\cI_*(\la, k)$ denote the sheaf ${}^p j_*(k[\dim \Gr^\la])$.  Given some other $\mu \in \La_+$, we have a natural morphism
\[ \cI_*(\la,k) \ast \cI_*(\mu,k) \to \cI_*(\la+\mu,k). \]

Let $V^\la = \bH^*(\Gr,\cI_*(\la,k))$ and $L^\la = F_\la(\cI_*(\la,k))$.  Then $L^\la$ is a $\check T_k$-invariant line in $V^\la$ such that
\[  \xymatrix{ L^\la \otimes L^\mu \ar[r] \ar[d] & L^{\la + \mu} \ar[d] \\ V^\la \otimes V^\mu \ar[r] & V^{\la+\mu} } \]
commutes for all $\la,\mu \in \La_+$.  Then $\check B_k$ is the subgroup of $Aut^\otimes(\bH^*)$ preserving the lines $L^\la$. In fact, we have a canonical basis vector $v_\la \in L^\la$ given by
\[  v_\la = [S_\la \cap \Gr^\la] \in F_\la(\cI_*(\la,k)) = H^{2 \dim(S_\la \cap \Gr^\la)}_c(S_\la \cap \Gr^\la)  \]
such that $v_\la \otimes v_\mu$ maps to $v_{\la+\mu}$. 

Recall that the weight lattice $X^*(\check T^\rho_k) $ is equal to the group of co-invariants $ X_*(T)_\Ga$ for the action of $\Ga$ on $X_*(T)$.  We will denote the natural projection as follows:
\[ X_*(T) \to X_*(T)_\Ga \]
\[ \la \mapsto \ov{\la} \]
Given a dominant $\ov{\la} \in X_*(T)_\Ga$, let
\[ V^{\ov \la} = \bigoplus_{\al \in \ov{\la} }\bH^*(\Gr, \cI_*(\al,K))^\Ga \]
and let
\[ L^{\ov{\la}} \su V^{\ov{\la}} \]
be the line generated by 
\[ v_{\ov{\la}} = \sum_{\al \in \ov{\la} } v_\al \]
Since 
\[ v_{\ov{\la}} \otimes v_{\ov{\mu}} \mapsto \sum_{\al \in \ov{\la}} \sum_{\be \in \ov{\mu}} v_{\al + \be}  = v_{\ov{\la + \mu}}\]
we conclude that there is a canonical Borel subgroup $\check T_ k \su \check B_k \su \check G_k$. \end{proof}

\se{Inner forms of $\check G_k$}\label{inner}

Let $K/k$ be a finite Galois extension with Galois group $\Ga$. Recall from example \ref{vect} that $Vect_K$ has a natural descent datum $(\ga_a, \id)$ such that $Vect_K^\Ga \simeq Vect_k$.  It turns out that we can twist this descent datum to produce other $K/k$-forms of $Vect_k$.  Let 
\[ \de: H^2(\Ga, K^\times) \simeq Br(K/k) \]
denote the natural isomorphism.

\begin{prop}\label{Brauer}
Let $[\zeta] \in H^2(\Ga, K^\times)$.  Then $(\ga_a, \zeta(a,b))$ defines a descent datum on $Vect_K$ such that $Vect_K^\Ga \simeq Mod_A$, where $A$ is a central simple $k$-algebra such that $[A] = \de[\zeta]$ and $Mod_A$ is the category of finitely generated right $A$-modules.
\end{prop}

\begin{proof}
We may assume that $\zeta$ is a normalized $2$-cocycle. Define the descent datum corresponding to a normalized $\zeta \in Z^2(\Ga,K^\times)$ as follows: to each $a\in \Ga$ we associate the same $\ga_a$ defined in example \ref{vect}.  Given $V \in Vect_K$, we define the natural isomorphism $\ga_{ab}V \simeq \ga_a \ga_b V$ to be the following map:
\[ \zeta(a,b)\cdot \id_V : \ga_{ab}V \to \ga_a \ga_b V \] In order for this to be a descent datum, the following diagram must commute for all $a,b,c \in \Ga$:
\[ \xymatrix{  \ga_{a b c}(V) \ar[rr]^{\zeta(a,bc)} \ar[d]_{\zeta(ab,c)} & & \ga_a\ga_{bc}(V) \ar[d]^{a \zeta(b,c)} \\ \ga_{ab}\ga_c(V) \ar[rr]_{\zeta(a,b)} & & \ga_a\ga_b\ga_c(V)} \]
We see that this is exactly the condition for $\zeta$ to be a 2-cocycle.

Now suppose that this $\Ga$-action is equivalent to $(\ga_a,\id)$.  This means that we have an equivalence $F:Vect_K \simeq Vect_K$ and natural isomorphisms $\psi_a: \ga_a F  \to F \ga_a$ such that the following diagram commutes:
\[ \xymatrix{  \ga_{ab} F(V) \ar[r]^{\id_{F(V)}}  \ar[d]_{\psi_{ab}} & \ga_a \ga_b F(V) \ar[r]^{\ga_a(\psi_b)} & \ga_a F(\ga_b V)\ar[d]^{\psi_a} \\ F(\ga_{ab} V) \ar[rr]_{F(\zeta(a,b))}&  & F(\ga_a\ga_b V) } \]
In other words,
\[ \psi_a \cdot a (\psi_b) = f(a,b) \cdot \psi(ab) \]
This shows that $\zeta$ comes from a map $\psi: \Ga \to K^\times$, hence is a $2$-coboundary.  Therefore, $H^2(\Ga, K^\times)$ is equal to the set of equivalence classes of descent data on $Vect_K$. 

Finally, recall the construction of an algebra $A$ in the Brauer class $\de[\zeta]$:
\[  A = \bigoplus_{a \in \Ga} K \cdot e_a \]
where the multiplication of the basis vectors is given by $e_a \cdot e_b = \zeta(a,b) e_{ab}$.  Given an object $(V,a_V)$ in $Vect_K^\Ga$, we see that
\[  e_a \mapsto a_V \]
is a $k$-algebra homomorphism from $A^{op} \to End_k(V)$.  This construction gives an equivalence $Vect_K^\Ga \simeq Mod_A$.
\end{proof}  

\bigskip

\begin{lem}\label{modules}
Let $K/k$ be a finite Galois extension with Galois group $\Ga$ and let $\cC$ be a $k$-linear abelian category.  Suppose we have central simple $k$-algebras $A$ and $B$ split by $K$, such that $[A] = [B] \in Br(K/k)$.  Then there is an equivalence of categories
\[ i_{A,B}: \cC_A \simeq \cC_B \]
where $\cC_A$ denotes the category of right $A$-modules in $\cC$.  If $C$ is another such algebra, then there is a natural isomorphism 
$\th(A,B,C) : i_{A,C} \simeq i_{B,C} i_{A,B} $.  Finally, if $D$ is another such algebra, then the following diagram commutes:
\[ \xymatrix{ i_{A,D} \ar[rrr]^{\th(A,C,D)} \ar[d]_{\th(A,B,D)} && & i_{C,D} i_{A,C}\ar[d]^{i_{C,D}(\th(A,B,C))} \\ i_{B,D} i_{A,B} \ar[rrr]_{\th(B,C,D)_{i_{A,B}(\un{\;\;})}} & && i_{C,D} i_{B,C} i_{A,B} } \]
\end{lem}

\begin{proof}
There is a natural descent datum $(\ga_a, \id)$ on $\cC_K$, the category of $K$-modules in $\cC$, such that 
\[ (\cC_K)^{(\ga_a,\;\id)} \simeq \cC\]
Let $\zeta_A \in Z^2(\Ga,K^\times)$ such that $\de[\zeta_A] = [A]$.  The proof of Proposition \ref{Brauer} implies that the descent datum $(\ga_a, \zeta(a,b))$ on $\cC_K$ has the following property:
\[ (\cC_K)^{(\ga_a,\;\zeta_A(a,b))} \simeq \cC_A \]
Therefore, for any central simple $k$-algebra $B$ represented by $\zeta_B \in Z^2(\Ga,K^\times)$ such that $[A] = [B] \in Br(K/k)$, we have an equivalence 
\[ i_{A,B} : \cC_A \simeq (\cC_K)^{(\ga_a, \; \zeta_A(a,b))} \simeq (\cC_K)^{(\ga_a, \; \zeta_B(a,b))} \simeq \cC_B \]
induced by the 2-coboundary $b_{A,B} = \zeta_A - \zeta_B \in B^2(\Ga, K^\times)$.  Given another such algebra $C$ represented by $\zeta_C$, we have:
\[ b_{A,C} = \zeta_A - \zeta_C = (\zeta_A - \zeta_B) + (\zeta_B - \zeta_C) = b_{A,B} + b_{B,C} \]
which induces an isomorphism
\[ \th(A,B,C) : i_{A,C} \simeq i_{B,C} i_{A,B} \]
Therefore, given another such algebra $D$ represented by $\zeta_D$, we have that
\[ b_{A,D} = b_{A,C} + b_{C,D} = b_{A,B} + b_{B,D} = b_{A,B} + b_{B,C} + b_{C,D} \]
These equations immediately imply that the corresponding diagrams commute.
\end{proof}

Let $G$ be a complex connected reductive algebraic group, $k$ a field, and $\mu : \pi_1(G)\to Br(k)$ a group homomorphism.  Recall that $\pi_0(\Gr) = \pi_1(G)$, so we can write \[ \Gr = \bigsqcup_{\al \in \pi_1(G)} \Gr_\al .\]
Let $P(\Gr,\mu)$ denote the following category:
\[ P(\Gr, \mu) := \bigoplus_{\al \in \pi_1(G)}P_{G_\cO}(\Gr_\al,A_\al) \]
where $A_\al$ is a central simple $k$-algebra such that $[A_\al]  = \mu(\al) \in Br(K/k)$ for some finite Galois extension $K/k$, and $P_{G_\cO}(\Gr_\al, A_\al)$ denotes the $k$-linear category of right $A_\al$-modules in $P_{G_\cO}(\Gr_\al,K)$.  

Let us describe the tensor structure on $P(\Gr,\mu)$. Fix a representative $A_\al$ for $\mu(\al)\in Br(K/k)$, for each $\al \in \pi_1(G)$. Let $\cA \in P_{G_\cO}(\Gr_\al,A_\al)$ and $\cB \in P_{G_\cO}(\Gr_\be,A_\be)$ and let $\cA \ast \cB$ denote the convolution of the underlying perverse $K$-sheaves.  The result is a right $A_\al \otimes A_\be$-module in $P_{G_\cO}(\Gr_{\al + \be}, K)$.  Since $[A_\al \otimes A_\be] = [A_{\al + \be}]$, Lemma \ref{modules} implies that we have an equivalence 
\[i_{\al ,\be}: P_{G_\cO}(\Gr_{\al + \be}, A_\al \otimes A_\be) \simeq P_{G_\cO}(\Gr_{\al + \be}, A_{\al + \be})\]
We define the tensor structure on $P(\Gr,\mu)$ as follows:
\[ \star: P(\Gr,\mu) \times P(\Gr,\mu) \to P(\Gr,\mu) \]
\[ \cA_\al \star \cA_\be = i_{\al,\be}(\cA_\al \ast \cA_\be). \]
The fact that this tensor structure is well-defined follows from Lemma \ref{modules}.

\begin{thm} If $\tilde G_k$ is an inner form of $\check G_k$, then $Rep_k(\tilde G_k)$ is tensor equivalent to $P(\Gr,\mu)$ for some $\mu: \pi_1(G) \to Br(k)$.
\end{thm}

\begin{proof}
Let $K/k$ be a finite Galois extension such that $\check G_K \simeq \tilde G_K$. Let $\Ga$ be the Galois group of $K/k$. Recall that $\tilde G_k$ is represented by a one-cocyle $[c] \in H^1(\Ga, Int(\check G_K))$.  Let $[\zeta] \in H^2(\Ga, Z(\check G_K))$  denote the image of $[c]$ under the boundary map.  Let $\mu$ denote the image of $[\zeta]$ under the isomorphism from Proposition \ref{H2}:
\[ H^2(\Ga, Z(\check G_K)) \simeq Hom(\pi_1(G), Br(K/k))\] 
\[ [\zeta] \mapsto \mu \]
We can use the cocycle $[c] \in H^1(\Ga, Int(\check G_K))$ to build a tensor descent datum on $P_{G_\cO}(\Gr,K)$ such that
\[ P_{G_\cO}(\Gr,K)^\Ga \simeq Rep(\tilde G_k)\]
As in Remark \ref{sheaf}, this descent datum is given by 
\[ (c(a\inv)^* \tilde \ga_a,\;\; \xi(c(b\inv), c(a \inv)) ) \]
We now define an equivalent descent datum on $P_{G_\cO}(\Gr,K)$ as follows:
\[ (\tilde \ga_a, \;\; \zeta(a,b)\cdot \id ) \]
The fact that this is a tensor descent datum follows from the following enlargement of diagram \ref{Satake}:
\begin{equation}\label{Satake} \xymatrix{ P_{G_\cO}(\Gr,K) \ar[rr]^\simeq \ar[d]_{\bigoplus_{\nu} F_\nu} && Rep(\check G_K) \ar[d] \\ Vect_K(X_*(T)) \ar[rr]^\simeq  \ar[d] && Rep(\check T_K) \ar[d] \\  Vect_K(\pi_1(G)) \ar[rr]^\simeq&& Rep(Z(\check G_K)) } \end{equation}
which shows that elements of $Z(\check G_K)$ act by tensor automorphisms. We have used the identification $\pi_1(G) \simeq X_*(T)/ X_*(T_{sc})$, where $T_{sc}$ is the pre-image of $T$ under the natural morphism $G_{sc} \to G$.

The fact that these two descent data are equivalent follows from the fact that
\[ \zeta(a,b) = c(a) \cdot c(b) \cdot c(ab)\inv \]
In other words, a diagram of type \ref{eta} commutes, with $F = \id$ and $\eta(a) : c(a\inv)^* \tilde \ga_a \simeq \tilde \ga_a$ is induced by $c(a)$.  Consequently, we have an equivalence of tensor categories
\[ P_{G_\cO}(\Gr,K)^{(c(a\inv)^* \tilde \ga_a,\;\; \xi(c(b\inv), c(a \inv)) )} \simeq P_{G_\cO}(\Gr,K)^{(\tilde \ga_a, \;\; \zeta(a,b)\cdot \id )} \]
Let us denote the latter category as $P(\Gr,\zeta)$. We will describe the objects of $P(\Gr,\zeta)$ that are supported on each component $\Gr_\al \su \Gr$, for $\al \in \pi_1(G)$.  Recall that the cocycle $\zeta$ is a map
\[ \Ga \times \Ga \to Z(\check G_K) = Hom(\pi_1(G),K^\times) \]
so that, for a fixed $\al \in \pi_1(G)$, we have $\zeta(\un{\;\;},\un{\;\;})(\al) \in Z^2(\Ga,K^\times)$. Let $(\cF, a_{\cF}) \in P(\Gr, \zeta)$ such that $\cF$ is supported on $\Gr_\al$.  This means that the following diagram commutes:
\[  \xymatrix{ \cF \ar[rr]^{a_{\cF}} \ar[d]_{(ab)_\cF} && \tilde \ga_a \cF \ar[d]^{\tilde \ga_a(b_\cF)} \\  \tilde \ga_{ab} \cF \ar[rr]^{\zeta(a,b)(\al) \cdot \id_\cF} && \tilde \ga_a \tilde \ga_b \cF} \]
Following the proof of Proposition \ref{Brauer}, this implies that $(\cF, a_\cF)$ defines a right $A_\al$-module $(\cF,r)$ in $P_{G_\cO}(\Gr_\al,K)$, where
\[ r: A_\al^{op} \to End(\cF) \]
\[ e_a \mapsto a_\cF \]
As usual, we have made use of the notation
\[  A_\al = \bigoplus_{a \in \Ga} K \cdot e_a \]
where the multiplication of the basis vectors is given by $e_a \cdot e_b = \zeta(a,b)(\al) \cdot e_{ab}$.  Therefore, we have constructed an equivalence of $k$-linear abelian categories:
\[ P(\Gr, \zeta) \simeq P(\Gr,\mu). \]
A version of Lemma \ref{modules} for tensor categories implies that this identification is an equivalence of tensor categories.
\end{proof}

\begin{rmk}
Note that $P(\Gr,\mu)$ has a fibre functor with values in $ Vect_K $, defined as the composition
\[ P(\Gr, \mu) \simeq P_{G_\cO}(\Gr,K)^\Ga \to P_{G_\cO}(\Gr,K) \to Vect_K \]
Given a rigid abelian tensor category $\cC$ over $k$ with a $K$-valued fibre functor, we can define the stack 
\[ Fib(\cC): k\tn{-alg} \to Gpd \]
such that $Fib(\cC)(R)$ is the groupoid of all $R$-valued fibre functors on $\cC$.  In fact, $Fib(\cC)$ is a gerbe (see \cite{DM}, 3.6 - 3.10). In our situation, this theorem implies that that $Fib(P(\Gr,\mu))$ is equivalent to the gerbe associated to the inner form $\tilde G_k$.
\end{rmk}

\end{document}